\documentclass[11pt]{article}%
\usepackage{amsfonts}
\usepackage{graphicx}
\usepackage{amsmath}
\usepackage{amssymb}%
\providecommand{\U}[1]{\protect\rule{.1in}{.1in}}
\newtheorem{theorem}{Theorem}

\newtheorem{conjecture}[theorem]{Conjecture}
\newtheorem{corollary}[theorem]{Corollary}

\newtheorem{lemma}[theorem]{Lemma}

\newtheorem{proposition}[theorem]{Proposition}
\newtheorem{remark}[theorem]{Remark}

\begin{document}

\begin{center}
{\Large On the functional CLT for stationary Markov Chains started at a point}

\bigskip

Dedicated to the memory of Mikhail Gordin\bigskip

David Barrera, Costel Peligrad and Magda Peligrad

\bigskip

\bigskip

\end{center}

Department of Mathematical Sciences, University of Cincinnati, PO Box 210025,
Cincinnati, Oh 45221-0025, USA.

E-mail: barrerjd@mail.uc.edu; peligrc@ucmail.uc.edu;
peligrm@ucmail.uc.edu\bigskip



\begin{center}
Abstract
\end{center}

We present a general functional central limit theorem started at a point also
known under the name of quenched. As a consequence, we point out several new
classes of stationary processes, defined via projection conditions, which
satisfy this type of asymptotic result. One of the theorems shows that if a
Markov chain is stationary ergodic and reversible, this result holds for
bounded additive functionals of the chain which have a martingale coboundary
in $\mathbb{L}_{1}$ representation. Our results are also well adapted for
strongly mixing sequences providing for this case an alternative, shorter
approach to some recent results in the literature.

\section{Introduction and results}

In this paper we address the question of the validity of functional limit
theorem for processes started at a point for almost all starting points. These
types of results are also known under the name of quenched limit theorems or
almost sure conditional invariance principles. The quenched functional CLT is
more general than the usual one and it is very important for analyzing random
processes in random environment, Markov chain Monte Carlo procedures and the
discrete Fourier transform (see Rassoul-Agha and Sepp\"{a}l\"{a}inen 2007,
2008, Barrera and Peligrad, 2016). On the other hand there are numerous
examples of processes satisfying the functional CLT\ but failing to satisfy
the quenched CLT. Some examples were constructed by Voln\'{y} and Woodroofe
(2010) and for the discrete Fourier transforms by Barrera (2015). This is the
reason why it is desirable to point out classes of processes satisfying a
quenched CLT. Special attention will be devoted to reversible Markov chains
and several open problems will be pointed out. Reversible Markov chains have
applications to statistical mechanics and to Metropolis Hastings algorithms
used in Monte Carlo simulations. The methods of proof we used are based on
martingale techniques combined with results from ergodic theory.

The field of limit theorems for stationary stochastic processes is closely
related to Markov operators and dynamical systems. All the results for
stationary sequences can be translated in the language of Markov operators and
vice-versa. In this paper we shall mainly use the Markov operator language and
also indicate the connection with stationary processes.

We assume that $(\xi_{n})_{n\in\mathbb{Z}}$ is a stationary Markov chain
defined on a probability space $(\Omega,\mathcal{F},\mathbb{P})$ with values
in a measurable state space $(S,\mathcal{A}),$ with marginal distribution
$\pi(A)=\mathbb{P}(\xi_{0}\in A)$ and regular conditional distribution for
$\xi_{1}$ given $\xi_{0},$ denoted by $Q(x,A)=\mathbb{P}(\xi_{1}\in A|\xi
_{0}=x)$. Let $Q$ also denote the Markov operator {acting via $(Qf)(x)=\int
_{S}f(s)Q(x,ds).$ Next, for }$p\geq1,$ {let $\mathbb{L}_{p}^{0}(\pi)$ be the
set of measurable functions on $S$ such that $\int|f|^{p}d\pi<\infty$ and
$\int fd\pi=0.$ For some function }${f}\in${$\mathbb{L}_{2}^{0}(\pi)$, let}%
\begin{equation}
{X_{i}=f(\xi_{i}),\ S_{n}=S_{n}(f)=\sum\limits_{i=1}^{n}X_{i}}. \label{defX}%
\end{equation}
{\ Denote by $\mathcal{F}_{k}$ the $\sigma$--field generated by $\xi_{i}$ with
$i\leq k$. }For any integrable random variable $X$ we denote by $\mathbb{E}%
_{k}(X)=\mathbb{E}(X|\mathcal{F}_{k})$ the conditional expectation of $X$
given $\mathcal{F}_{k}.$ With this notation, $\mathbb{E}_{0}(X_{1}%
)=(Qf)(\xi_{0})=\mathbb{E}(X_{1}|\xi_{0}).$ We denote by ${{||X||}_{p}}$ the
norm in {$\mathbb{L}_{p}=\mathbb{L}_{p}$}$(\Omega,\mathcal{F},\mathbb{P}).$
The integral on the space $(S,\mathcal{A},\pi)$ will be denoted by
$\mathbb{E}_{\pi}.$ So, $\mathbb{E}{f(\xi_{0})=}\mathbb{E}_{\pi}f.$

The Markov chain is usually constructed in a canonical way on $\Omega
=S^{\infty}$ endowed with sigma algebra $\mathcal{A}^{\infty},$ and
$\mathcal{\xi}_{n}$ is the $n^{th}$ projection on $S$. The shift
$T:\Omega\rightarrow\Omega$ is defined by $\mathcal{\xi}_{n}(T\omega
)=\mathcal{\xi}_{n+1}(\omega)$ for every integer $n.$

For any probability measure $\upsilon$ on $\mathcal{A}$ the law of $(\xi
_{n})_{n\in\mathbb{Z}}$ with transition operator $Q$ and initial distribution
$\upsilon$ is the probability measure $\mathbb{P}^{\upsilon}$ on $(S^{\infty
},\mathcal{A}^{\infty})$ such that%
\[
\mathbb{P}^{\upsilon}(\xi_{n+1}\in A|\xi_{n}=x)=Q(x,A)\text{ and }%
\mathbb{P}^{\upsilon}(\xi_{0}\in A)=\upsilon(A).
\]
For $\upsilon=\pi$ we denote $\mathbb{P}=\mathbb{P}^{\pi}.$ For $\upsilon
=\delta_{x},$ the Dirac measure, we denote by $\mathbb{P}^{x}$ and
$\mathbb{E}^{x}$ the probability and conditional expectation for the process
started at $x$. Note that for each $x$ fixed $\mathbb{P}^{x}(\cdot)$ is a
measure on {$\mathcal{F}^{\infty},$ the sigma algebra generated by }$\cup_{k}%
${$\mathcal{F}_{k}.$ Also }%
\begin{equation}
\mathbb{P}(A)=%
{\displaystyle\int}
\mathbb{P}^{x}(A)\pi(dx). \label{erg dec}%
\end{equation}

We mention that any stationary sequence $(Y_{k})_{k\in\mathbb{Z}}$ can be
viewed as a function of a Markov process $\xi_{k}=(Y_{j};j\leq k)$ with the
function $g(\xi_{k})=Y_{k}$. Therefore the theory of stationary processes can
be imbedded in the theory of Markov chains. So, our results apply to any
stationary process with the corresponding interpretation. In the context of a
stationary process, a fixed starting point for a corresponding Markov chain
means a fixed past trajectory for $k\leq0$.

All along the paper we shall assume that the Markov chain is ergodic.

Below, we denote by $\Rightarrow$ the convergence in distribution. By $[x]$ we
denote the integer part of $x.$

For a Markov chain, by the quenched CLT\ (or CLT started at a point) we shall
understand the following convergence: there is a positive constant $\sigma
\in\lbrack0,\infty)$ and a set $S^{\prime}\subset S$ with $\pi(S^{\prime})=1$
such that$\ $for $x\in S^{\prime}$ we have%
\begin{equation}
\frac{S_{n}}{\sqrt{n}}\Rightarrow\sigma N(0,1)\text{ under }\mathbb{P}^{x},
\label{QCLT}%
\end{equation}
and by the quenched functional CLT\ (which is the same as functional
CLT\ started at a point): there is a set $S^{\prime}\subset S$ with
$\pi(S^{\prime})=1$ such that$\ $for $x\in S^{\prime}$
\begin{equation}
\frac{S_{[nt]}}{\sqrt{n}}\Rightarrow\sigma W(t)\text{ under }\mathbb{P}^{x},
\label{FQCLT}%
\end{equation}
where $W(t)$ denotes the standard Brownian motion and the convergence in
distribution is on $D(0,1)$, the space of functions continuous at the right
with limits at the left, endowed with the Skorohod topology.

An important class satisfying quenched functional CLT is the stationary and
ergodic martingale differences, as seen in Derriennic and Lin (2001, 2003). A
natural method to prove these types of results for other classes of processes
is to use martingale approximations. This method was initiated by Gordin (1969).

One of the first results of this type is due to Gordin (published in Ch.4
Section 8 in Borodin and Ibragimov, 1994), who proved the quenched CLT\ for
Markov chains with normal operator ($QQ^{\ast}=Q^{\ast}Q$)$,$ $f\in
\mathbb{L}_{2}^{0},$ under the condition $f\in(I-Q)\mathbb{L}_{2}(\pi)$. If
the Markov chain is irreducible and aperiodic, then the quenched CLT holds
under the condition $\sum_{j=0}^{n}\mathbb{E}_{\pi}(fQ^{j}f)$ is convergent
(Chen, 1999). Without assuming irreducibility conditions, various papers point
out rates for convergence to $0$ of $||\sum_{j=0}^{n}Q^{j}f||_{2}/n$ needed
for the quenched results. Among them, we mention papers by Derriennic and Lin
(2001, 2003), Wu and Woodroofe (2004), Cuny (2011), Merlev\`{e}de et al.
(2011), Cuny and Peligrad (2012), Cuny and Merlev\`{e}de (2014), Cuny and
Voln\'{y} (2013), Voln\'{y} and Woodroofe (2014). Recently, Dedecker et al.
(2014) showed that the condition $\sum_{j=0}^{\infty}\mathbb{E}_{\pi}%
|fQ^{j}f|<$ $\infty\ $leads to the quenched invariance principle.

Our study is motivated by the class considered by Gordin. What can one say
about $f\in(I-Q)\mathbb{L}_{p}(\pi)$ with $1\leq p<2?$ From the paper by
Voln\'{y} and Woodroofe (2014) we know that there are examples of functions,
$f\in\lbrack(I-Q)\mathbb{L}_{1}(\pi)]\cap$ $\mathbb{L}_{2}^{0}(\pi)$ such that
$S_{n}/\sqrt{n}$ satisfies the CLT, but fails to satisfy the quenched CLT.

One of our results shows that for functions of reversible Markov chains one
can assume that $f\in\lbrack(I-Q)\mathbb{L}_{q}(\pi)]\cap\mathbb{L}_{p}%
^{0}(\pi),$ with $q\in\lbrack1,2],$ $1/p+1/q=1,$ for concluding that the
quenched functional CLT holds. This result follows from several general
preliminary results that have interest in themselves. They specify sufficient
conditions for the validity of the quenched CLT\ and the quenched functional CLT.

Denote
\begin{equation}
f_{m}=\frac{1}{m}\mathbb{(}Q+...+Q^{m})f \label{deffm}%
\end{equation}
and%
\begin{equation}
\bar{R}_{k}^{m}=\sum\nolimits_{j=1}^{k}f_{m}(\xi_{j})\text{.} \label{defRbar}%
\end{equation}

\begin{theorem}
$\ $\label{quenched} Let $(X_{n})_{n\in\mathbb{Z}}$ be a stationary sequence
of random variables defined by (\ref{defX})\ and define $(\bar{R}_{k}%
^{m})_{m\geq1,k\geq1}$ by (\ref{defRbar}). Assume that
\begin{equation}
\lim_{m}\lim\sup_{n}\mathbb{P}^{x}(\frac{|\bar{R}_{n}^{m}|}{\sqrt{n}%
}>\varepsilon)=0\text{ \ }\pi-\text{a.s.} \label{negl1Th}%
\end{equation}
Then the quenched CLT in (\ref{QCLT}) holds.
\end{theorem}

\begin{theorem}
\label{quenched-IP}Assume that $(X_{n})_{n\in\mathbb{Z}}$ and $(\bar{R}%
_{k}^{m})_{m\geq1,k\geq1}$ are as in Theorem \ref{quenched} and
\begin{equation}
\lim_{m}\lim\sup_{n}\mathbb{P}^{x}(\max_{1\leq j\leq n}\frac{|\bar{R}_{j}%
^{m}|}{\sqrt{n}}>\varepsilon)=0\text{ \ }\pi-\text{a.s.} \label{neglIPTh}%
\end{equation}
then the quenched functional CLT\ in (\ref{FQCLT})\ holds.
\end{theorem}

For ${f}\in${$\mathbb{L}_{1}^{0}(\pi)$ }denote by
\begin{equation}
g_{f}=\sup_{n\geq0}|\sum_{j=0}^{n}Q^{j}f|. \label{defg}%
\end{equation}
Based on Theorem \ref{quenched-IP} we shall establish the following theorem:

\begin{theorem}
\label{Th-quenched} Let $(X_{n})_{n\in\mathbb{Z}}$ be defined by (\ref{defX}),
$f_{m}$ by (\ref{deffm}) and $g_{f}$ by (\ref{defg}). Assume the following
condition is satisfied:
\begin{equation}
(f_{m}g_{f})_{m\geq1}\text{ is uniformly integrable.} \label{Qnew}%
\end{equation}
Then the quenched functional CLT in (\ref{FQCLT}) holds.
\end{theorem}

From the proof of Theorem \ref{Th-quenched} we easily deduce several
corollaries. The first corollary is well adapted for strongly mixing sequences:

\begin{corollary}
\label{corstrong}Assume%
\begin{equation}
\lim_{m\rightarrow\infty}\sum_{j=1}^{\infty}\mathbb{E}_{\pi}|(Q^{m}%
f)(Q^{j}f)|=0. \label{strong}%
\end{equation}
Then the quenched functional CLT\ in (\ref{FQCLT}) holds.
\end{corollary}

\begin{remark}
\label{strong mixing}Condition (\ref{strong})\ can be verified in terms of
strong mixing coefficients. Practically, we deduce that any strongly mixing
sequence satisfying the CLT also satisfies the quenched functional CLT.
Therefore our approach also provides a shorter, alternative proof of Corollary
3.5 in Dedecker et al. (2014). The proof of this remark is postponed to the
end of the paper.
\end{remark}

Also, as an application to the proof of Theorem \ref{Th-quenched} we obtain
the following:

\begin{corollary}
\label{cor-quenched}Let $(X_{n})_{n\in\mathbb{Z}}$, $f_{m}$, and $g_{f}$
defined as in Theorem \ref{Th-quenched}. Assume $f\in\mathbb{L}_{p}^{0}(\pi)$
and $g_{f}\in\mathbb{L}_{q}(\pi)$ with $p\in\lbrack2,\infty],$ $1/p+1/q=1.$
Then the quenched functional CLT holds.
\end{corollary}

We say that a Markov chain is reversible if $Q$ is self-adjoint; equivalently
$(X_{0},X_{1})$ and $(X_{1},X_{0})$ are identically distributed. If the Markov
chain is reversible then the following corollary holds.

\begin{corollary}
\label{rev}Assume the Markov chain is reversible and
\begin{equation}
f\in\lbrack(I-Q)\mathbb{L}_{q}(\pi)]\cap\mathbb{L}_{p}^{0}(\pi).
\label{coboundary}%
\end{equation}
for $p\in\lbrack2,\infty),$ $1/p+1/q=1.$ Then the quenched functional CLT holds.
\end{corollary}

Let us mention that the class we consider here is of independent interest when
compared to the projective condition used in Dedecker et al. (2014), namely
$\sum_{j=0}^{\infty}\mathbb{E}|X_{0}E(X_{j}|\mathcal{F}_{0})|<\infty$. For
instance there are examples which satisfy the conditions of Corollary
\ref{cor-quenched} without satisfying the condition from Dedecker et al. (2014).

\begin{remark}
\label{example}There is a stationary and ergodic process of bounded random
variables $(X_{k})_{k\in Z}$ adapted to a filtration $(\mathcal{F}_{k})_{k\in
Z}$, such that $\sup_{n\geq0}|\sum_{j=0}^{n}\mathbb{E}(X_{j}|\mathcal{F}%
_{0})|\in\mathbb{L}_{1}$ and$\ \sum_{j=0}^{\infty}\mathbb{E}|X_{0}%
\mathbb{E}(X_{j}|\mathcal{F}_{0})|=\infty.$
\end{remark}

We end this section by mentioning two conjectures which deserve further
investigation. The results in the paper by Dedecker et al. (2014) and the
results in this paper suggest the following conjecture, which is a quenched
form of the functional CLT in Dedecker and Rio (2000).

\begin{conjecture}
In the context of Theorem \ref{Th-quenched} assume
\begin{equation}
|f\sum\nolimits_{j=0}^{n}Q^{j}f|\text{ is convergent in }\mathbb{L}_{1}(\pi).
\label{DR}%
\end{equation}
Then the quenched functional CLT holds.
\end{conjecture}

For reversible Markov chains we would like to mention the Kipnis and Varadhan
(1986) conjecture, asking if their functional CLT is quenched. This conjecture
is still unsolved.

\begin{conjecture}
In the context of Corollary \ref{rev} assume
\begin{equation}
\mathbb{E}_{\pi}\mathbb{(}f\sum\nolimits_{j=0}^{n}Q^{j}f)\text{ is
convergent.} \label{KV}%
\end{equation}
Then the quenched functional CLT holds.
\end{conjecture}

Steps towards clarifying this conjecture are contained in the papers by
Derriennic and Lin (2001) and Cuny and Peligrad (2012).

\section{Preliminary considerations}

The method we shall use in our proofs\ is based on a martingale approximation
depending on a certain parameter which is fixed at the beginning and after
that we let it grow to $\infty$. To deal with this parameter, we start by
pointing out several preliminary considerations for convergence in
distribution. From Theorem 3.2 in Billingsley (1999), it is well-known the
following result:

\begin{lemma}
\label{Billingsley}Assume that the elements $(X_{n,m},X_{n})$ are defined on
the same probability space with values in $S\times S,$ where $S$ is a metric
space. Assume that
\[
X_{n,m}\Rightarrow_{n}Y_{m}\Rightarrow_{m}X
\]
and
\begin{equation}
\lim_{m}\lim\sup_{n}\mathbb{P}(d(X_{n,m},X_{n})\geq\varepsilon)=0.
\label{negl1}%
\end{equation}
Then%
\[
X_{n}\Rightarrow X.
\]

\end{lemma}

If the metric space is separable and complete, then one does not have to
assume $Y_{m}\Rightarrow X$. This result (see for instance Theorem 2 in
Dehling et al., 2009) is given in the following lemma where the variables are
denoted as in Lemma \ref{Billingsley}.

\begin{lemma}
\label{Dehling}Assume the metric space $S$ is separable and complete. Assume
that for every $m$
\[
X_{n,m}\Rightarrow Y_{m}\text{ as }n\rightarrow\infty
\]
and condition (\ref{negl1}) is satisfied. Then there is a $S$-valued random
variable $X$ such that
\[
Y_{n}\Rightarrow X\text{ and }X_{n}\Rightarrow X\text{ as }n\rightarrow
\infty.
\]

\end{lemma}

These considerations suggest that the conditions of Lemma \ref{Billingsley}
are too strong. Indeed, we can formulate the following lemma.

\begin{lemma}
\label{gen-bill}In Lemma \ref{Billingsley} condition \ref{negl1} can be
replaced by
\begin{equation}
\lim\inf_{m}\lim\sup_{n}\mathbb{P}(d(X_{n,m},X_{n})\geq\varepsilon)=0.
\label{negl2}%
\end{equation}

\end{lemma}

\textbf{Proof of Lemma \ref{gen-bill}}. Let $F$ be a closed set. Define
$F_{\varepsilon}=\{x:d(x,F)\leq\varepsilon\}$. Then, by Portmanteau Theorem
(Theorem 2.1 in Billingsley 1999),
\[
\lim\sup_{n}\mathbb{P}(X_{n,m}\in F_{\varepsilon})\leq\mathbb{P}(Y_{m}\in
F_{\varepsilon}).
\]
Since
\[
\mathbb{P}(X_{n}\in F)\leq\mathbb{P}(X_{n,m}\in F_{\varepsilon})+\mathbb{P}%
(d(X_{n,m},X_{n})\geq\varepsilon),
\]
by combining these results, we deduce that%
\begin{align*}
\lim\sup_{n}\mathbb{P}(X_{n}  &  \in F)\leq\lim\sup_{n}\mathbb{P}(X_{n,m}\in
F_{\varepsilon})+\lim\sup_{n}\mathbb{P}(d(X_{n,m},X_{n})\geq\varepsilon)\\
&  \leq\mathbb{P}(Y_{m}\in F_{\varepsilon})+\lim\sup_{n}\mathbb{P}%
(d(X_{n,m},X_{n})\geq\varepsilon).
\end{align*}
Therefore taking the limit inferior when $m\rightarrow\infty$ we obtain by
(\ref{negl2}) and Portmanteau Theorem that
\begin{align*}
\lim\sup_{n}\mathbb{P}(X_{n}  &  \in F)\leq\lim\inf_{m}[\mathbb{P}(Y_{m}\in
F_{\varepsilon})+\lim\sup_{n}\mathbb{P}(d(X_{n,m},X_{n})\geq\varepsilon
)]\leq\\
\lim\sup_{m}\mathbb{P}(Y_{m}  &  \in F_{\varepsilon})+\lim\inf_{m}\lim\sup
_{n}\mathbb{P}(d(X_{n,m},X_{n})\geq\varepsilon)\leq\mathbb{P}(X\in
F_{\varepsilon}).
\end{align*}
Now we take a sequence $F_{\varepsilon}\downarrow F$ as $\varepsilon
\downarrow0,$ the result follows by applying again the Portmanteau Theorem.
$\square$

\bigskip

One of the difficulties in proving quenched results is the fact that, under
$\mathbb{P}^{x},$ the Markov chain is no longer strictly stationary. Since we
are interested in proving quenched results which are almost sure results, and
also the quenched functional form of the CLT, we need to use maximal
inequalities. There are not too many maximal inequalities available in the
nonstationary context. A useful maximal inequality is an easy consequence of
inequality (3.9) given in the book by Rio (2000), (see also Dedecker and Rio, 2000).

\begin{lemma}
\label{Lemma rio}Assume that $(X_{k})$ is a sequence of real valued centered
random variables in $\mathbb{L}_{2}(\Omega,\mathcal{K},\mathbb{P}),$ adapted
to an increasing filtration of sub-sigma fields of $\mathcal{K}$,
$(\mathcal{F}_{n}).$ Then
\[
\mathbb{E}(\max_{1\leq k\leq n}S_{k}^{2})\leq8%
{\displaystyle\sum_{k=1}^{n}}
\mathbb{E}(X_{k}^{2})+16%
{\displaystyle\sum\limits_{k=1}^{n}}
\mathbb{E}|X_{k}\mathbb{E}(S_{n}-S_{k}|\mathcal{F}_{k})|.
\]

\end{lemma}

One of the basic results used in our proofs is the functional CLT for
martingale in the following form:

\begin{theorem}
\label{Mart QFCLT}Assume that $(D_{n})$ is a sequence of martingale
differences on a probability space $(\Omega,\mathcal{K},\mathbb{P})$ adapted
to an increasing filtration of sub-sigma fields of $\mathcal{K}$,
$(\mathcal{F}_{n}).$ Assume that the following two conditions hold%
\begin{equation}
(\frac{1}{\sqrt{n}}\max_{0\leq k\leq n}|D_{k}|)_{n\geq1}\text{ is uniformly
integrable} \label{Mart cond 1}%
\end{equation}
and for each $t,$ $0\leq t\leq1$%
\begin{equation}
\frac{1}{n}%
{\displaystyle\sum_{k=0}^{[nt]}}
D_{k}^{2}\rightarrow t\sigma^{2}\text{ in probability.} \label{mart cond 2}%
\end{equation}
Then
\[
\frac{%
{\displaystyle\sum_{k=0}^{[nt]}}
D_{k}}{\sqrt{n}}\Rightarrow|\sigma|W(t).
\]

\end{theorem}

This theorem follows from Theorem 2.3 in Gaenssler and Haeusler\ (1986)
combined with the commentaries on pages 316-317 of this paper. Indeed,
according to the sequence of implications on page 316 of this book, the
conditions (A$_{a}$) and (R$_{a,t}$) of their Theorem 2.3 are verified under
(\ref{mart cond 2}) and
\begin{equation}
\frac{1}{\sqrt{n}}\max_{0\leq k\leq n}|D_{k}|\rightarrow0\text{ in }%
\mathbb{L}_{1}. \label{mart l1 cond}%
\end{equation}
Then, by arguments on page 317 both conditions (\ref{Mart cond 1}) and
(\ref{mart cond 2}) imply condition (\ref{mart l1 cond}).

\section{Proofs}

\textbf{Proof of Theorems \ref{quenched} and \ref{quenched-IP}.}

We start with a martingale construction. The construction of the martingale
decomposition is inspired by works of Gordin (1969), Heyde (1974),
Gordin-Lifshitz (1981); see also Theorem 8.1 in Borodin and Ibragimov (1994),
and Kipnis and Varadhan (1986) and Maxwell and Woodroofe (2000). The form we
use here was initiated by Wu and Woodroofe (2004), and further exploited by
Zhao and Woodroofe (2008), Peligrad (2010), Gordin and Peligrad (2011) among
others. We briefly give it here for completeness.

We introduce a parameter, an integer $m\geq1$ (kept fixed for the moment), and
introduce the functions%
\begin{equation}
v_{k}=(I+Q+...+Q^{k-1})f. \label{def f}%
\end{equation}
Define the stationary sequence of random variables:%
\[
\theta_{0}^{m}=\frac{1}{m}\sum_{k=1}^{m}v_{k}(\xi_{0}),\text{ }\theta_{k}%
^{m}=\theta_{0}^{m}\circ T^{k}\text{.}%
\]
Denote by
\begin{equation}
D_{k}^{m}=D_{k}^{m}(\xi_{k},\xi_{k+1})=\theta_{k+1}^{m}-\mathbb{E}_{k}%
(\theta_{k+1}^{m})\text{ ; }M_{n}^{m}=\sum_{k=1}^{n}D_{k}^{m}\text{.}
\label{defD}%
\end{equation}
Then, $(D_{k}^{m})_{k\in\mathbb{Z}}$ is a martingale difference sequence which
is stationary and ergodic and $(M_{n}^{m})_{n\geq0}$ is a martingale. So we
have%
\[
X_{k}=D_{k}^{m}+\theta_{k}^{m}-\theta_{k+1}^{m}+f_{m}(\xi_{k}),
\]
with $f_{m}$ defined by (\ref{deffm}). Therefore%
\begin{equation}
S_{k}=M_{k}^{m}+\theta_{1}^{m}-\theta_{k+1}^{m}+\overline{R}_{k}^{m},
\label{martdec}%
\end{equation}
where we implemented the notation%
\[
\bar{R}_{k}^{m}=\sum\nolimits_{j=1}^{k}f_{m}(\xi_{j}).
\]
With the notation
\begin{equation}
R_{k}^{m}=\theta_{1}^{m}-\theta_{k+1}^{m}+\bar{R}_{k}^{m}, \label{def Rest}%
\end{equation}
we have the following martingale decomposition
\begin{equation}
S_{k}=M_{k}^{m}+R_{k}^{m}\text{.} \label{martingale decomposition}%
\end{equation}
We shall prove now the quenched functional CLT for the martingale\textbf{
}$M_{n}^{m}$. We shall verify the conditions of the functional CLT given in
Theorem \ref{Mart QFCLT}.

We start by noticing that $(M_{n}^{m})_{n}$ is also a martingale under
$\mathbb{P}^{x}$ (since $\mathbb{E}^{x}(D_{k}^{m}|\mathcal{F}_{k-1}%
)=\mathbb{E}(D_{k}^{m}|\mathcal{F}_{k-1})$ by the fact that the Markov chain
has the same transitions under $\mathbb{P}$ and $\mathbb{P}^{x}).$ We verify
first condition (\ref{mart cond 2}). Since $M_{n}^{m}$ is a martingale with
stationary and ergodic increments, by Birkhoff's ergodic theorem, for every
$0\leq t\leq1,$%
\[
\frac{1}{n}\sum_{k=1}^{[nt]}(D_{k}^{m})^{2}\rightarrow t\mathbb{E}(D_{0}%
^{m})^{2}\text{ \ }\mathbb{P}-a.s.
\]
and therefore for every $0\leq t\leq1$ and $\pi-$almost all $x$
\begin{equation}
\frac{1}{n}\sum_{k=1}^{[nt]}(D_{k}^{m})^{2}\rightarrow t\mathbb{E}(D_{0}%
^{m})^{2}\text{ \ }\mathbb{P}^{x}-a.s. \label{erg}%
\end{equation}
In order to verify (\ref{Mart cond 1}), for proving uniform integrability it
is enough to show that for $\pi-$almost all $x,$ for some constant $C_{x}$ we
have
\begin{equation}
\sup_{n}\frac{1}{n}\mathbb{E}^{x}(\max_{1\leq k\leq n}(D_{k}^{m})^{2})\leq
C_{x}. \label{mart cond 11}%
\end{equation}
Clearly%
\[
\frac{1}{n}\mathbb{E}^{x}(\max_{1\leq k\leq n}(D_{k}^{m})^{2})\leq\frac{1}%
{n}\sum_{k=1}^{n}\mathbb{E}^{x}(D_{k}^{m})^{2}.
\]
Note that $D_{0}^{m}=D_{0}^{m}(\xi_{1},\xi_{0})$ and then, denoting by
$h(y)=E((D_{0}^{m}(\xi_{1},\xi_{0}))^{2}|\xi_{0}=y),$ by the Markov property
it follows that $\mathbb{E}^{x}(D_{k}^{m})^{2}=Q^{k}h(x).$ By Hopf's ergodic
theorem for Markov operators (see Theorem 11.4 in Eisner et al. 2015) we
obtain
\[
\lim\sup_{n}\frac{1}{n}\mathbb{E}^{x}(\max_{1\leq k\leq n}(D_{k}^{m})^{2}%
)\leq\lim\sup_{n}\frac{1}{n}\sum_{k=1}^{n}Q^{k}h(x)=\mathbb{E}(D_{0}^{m}%
)^{2}\text{ \ }\mathbb{\pi}-a.s.
\]
and (\ref{mart cond 11}) follows.

By Theorem \ref{Mart QFCLT} it follows that for $\pi-$almost all $x$ we have
\begin{equation}
\frac{M_{[nt]}^{m}}{\sqrt{n}}\Rightarrow|\sigma_{m}|W(t)\text{ under
}\mathbb{P}^{x}, \label{FCLTmart}%
\end{equation}
where $W(t)$ is the standard Brownian motion and
\begin{equation}
\sigma_{m}^{2}=\mathbb{E}(D_{0}^{m})^{2}. \label{defsigmam}%
\end{equation}
By stationarity, by the fact that $\theta_{0}^{m}$ is in $\mathbb{L}_{2}$ we
have
\[
\frac{\max_{1\leq k\leq n}|\theta_{k}^{m}|}{\sqrt{n}}\rightarrow0\text{
\ }\mathbb{P}-a.s.
\]
To see it, just start from $%
{\displaystyle\sum\nolimits_{n}}
\mathbb{P}(|\theta_{0}^{m}|^{2}>\varepsilon n)<\infty$ and apply the
Borel-Cantelli lemma (see also page 171 in Borodin and Ibragimov, 1994).

Therefore, for $\pi-$almost all $x$
\begin{equation}
\frac{\max_{1\leq k\leq n}|\theta_{k}^{m}|}{\sqrt{n}}\rightarrow0\text{
\ }\mathbb{P}^{x}-a.s. \label{IPM2}%
\end{equation}
If we assume (\ref{negl1Th}) then clearly by (\ref{IPM2}) we obtain
\[
\lim_{m}\lim\sup_{n}\mathbb{P}^{x}(\frac{|S_{n}-M_{n}^{m}|}{\sqrt{n}%
}>\varepsilon)=0\text{ \ }\pi-a.s.
\]
Clearly (\ref{FCLTmart}) implies that for each $\ m\geq1$
\[
\frac{M_{n}^{m}}{\sqrt{n}}\Rightarrow|\sigma_{m}|Z\text{ under }\mathbb{P}^{x}%
\]
where $Z$ has a standard normal distribution. By applying Lemma \ref{Dehling}
we obtain that $|\sigma_{m}|Z$ converges in distribution to a random variable
$Y,$ which is also the limiting distribution of $S_{n}/\sqrt{n}$ under
$\mathbb{P}^{x}.$ Clearly $Y$ has a normal distribution with variance
$\sigma^{2}=\lim_{m}\sigma_{m}^{2},$ where $\sigma\in\lbrack0,\infty).$

Now, by taking into account (\ref{neglIPTh}), we have%
\[
\lim_{m}\lim\sup_{n}\mathbb{P}^{x}(\max_{1\leq j\leq n}\frac{|S_{n}-M_{n}%
^{m}|}{\sqrt{n}}>\varepsilon)=0\text{ \ }\pi-a.s.
\]
and, by Lemma \ref{Dehling}, as explained before, we get both that
$\mathbb{E}(D_{0}^{m})^{2}\rightarrow\mathbb{\sigma}^{2}$ and that the
quenched functional CLT holds with the limit $\sigma W(t)$. \ $\square$

\begin{remark}
We point out that the proof of Theorems \ref{quenched} and \ref{quenched-IP}
also indicates how to identify the constant $\sigma^{2}$ which appears in the
limit as
\[
\sigma^{2}=\lim_{m\rightarrow\infty}\lim_{n\rightarrow\infty}\mathbb{E}%
(D_{n}^{m})^{2},
\]
where $D_{n}^{m}$ was defined in (\ref{defD}).
\end{remark}

\begin{remark}
By Lemma \ref{gen-bill}, Theorems \ref{quenched} and \ref{quenched-IP}\textbf{
}also hold if we replace in conditions (\ref{negl1Th}) and (\ref{neglIPTh})
the limit when $m\rightarrow\infty$ by $\lim\inf_{m\rightarrow\infty}$ and we
add the condition
\begin{equation}
\mathbb{E}(D_{0}^{m})^{2}\text{ is convergent as }m\rightarrow\infty\text{.}
\label{cond heyde}%
\end{equation}
Condition (\ref{cond heyde}) is verified in many situations including classes
of normal and reversible Markov chains as shown by Gordin and Lifshitz (1981),
and Kipnis and Varadhan (1986) among others.
\end{remark}

We shall establish next a maximal inequality needed to verify condition
(\ref{neglIPTh}).

\begin{proposition}
\label{inequality} For any $h\in${$\mathbb{L}_{2}^{0}(\pi)$ such that
}$\mathbb{E}_{\pi}(|hg_{h})|)<\infty,$ we have the following maximal
inequality%
\begin{equation}
\lim\sup_{n}\frac{\mathbb{E}^{x}(\max_{1\leq k\leq n}S_{k}^{2}(h))}{n}%
\leq24\mathbb{E}_{\pi}(|hg_{h}|)\text{ \ }\pi-a.s. \label{ineq2}%
\end{equation}

\end{proposition}

\textbf{Proof. }We start by applying Rio's maximal inequality given in Lemma
\ref{Lemma rio} which implies%
\[
\mathbb{E}^{x}(\max_{1\leq k\leq n}S_{k}^{2}(h))\leq8%
{\displaystyle\sum\limits_{u=1}^{n}}
\mathbb{E}^{x}(h^{2}(\xi_{u}))+16%
{\displaystyle\sum\limits_{u=1}^{n-1}}
\mathbb{E}^{x}|h(\xi_{u})%
{\displaystyle\sum\limits_{k=1}^{n-u}}
Q^{k}h(\xi_{u})|.
\]
So%
\[
\mathbb{E}^{x}(\max_{1\leq k\leq n}S_{k}^{2}(h))\leq24\sum_{j=1}^{n}Q^{j}%
[\sup_{k\geq0}|\sum_{u=0}^{k}hQ^{u}h|](x).
\]
By the Hopf ergodic theorem for Markov operators
\[
\frac{1}{n}\sum_{j=1}^{n-1}Q^{j}[\sup_{k\geq0}|\sum_{u=0}^{k}hQ^{u}%
h|](x)\rightarrow\mathbb{E}_{\pi}\sup_{n\geq0}|\sum_{u=0}^{n}hQ^{u}h|\text{
\ }\pi-a.s.
\]
which leads by the previous considerations to (\ref{ineq2}) by the definition
of $g_{h}$. $\square$

\bigskip

\textbf{Proof of Theorem \ref{Th-quenched}.}

The proof consists in verifying condition (\ref{neglIPTh}) of Theorem
\ref{quenched-IP}.

We start by applying Proposition \ref{inequality} to $S_{k}(f_{m}),$ where
$f_{m}$ is defined by (\ref{deffm}). Note that $\bar{R}_{k}^{m}$ defined by
(\ref{defRbar}) is equal to $S_{k}(f_{m}).$\ For all $m$ fixed%
\begin{equation}
\lim\sup_{n}\frac{\mathbb{E}^{x}(\max_{1\leq k\leq n}(\bar{R}_{k}^{m})^{2}%
)}{n}\leq24\mathbb{E}_{\pi}[\sup_{k\geq0}|\sum_{j=0}^{k}f_{m}Q^{j}%
f_{m}|]\text{ \ }\pi-a.s. \label{ineq11}%
\end{equation}
Then, we have%
\begin{align*}
|\sum_{j=0}^{n}Q^{j}f_{m}|  &  =\frac{1}{m}|\sum_{j=0}^{n}%
{\displaystyle\sum\limits_{k=1}^{m}}
Q^{j+k}f|\leq\frac{1}{m}%
{\displaystyle\sum\limits_{k=1}^{m}}
|\sum_{j=k}^{n+k}Q^{j}f|\\
&  \leq2\sup_{n}|\sum_{j=0}^{n}Q^{j}f|\leq2g_{f},
\end{align*}
which, combined with (\ref{ineq11}), leads to
\[
\lim\sup_{n}\frac{\mathbb{E}^{x}(\max_{1\leq k\leq n}(\bar{R}_{k}^{m})^{2}%
)}{n}\leq48\mathbb{E}_{\pi}(|f_{m}g_{f}|)\text{ \ }\pi-a.s.
\]
Clearly, by using this last inequality, in order to prove (\ref{neglIPTh}), it
remains to show%
\begin{equation}
\mathbb{E}_{\pi}(|f_{m}g_{f}|)\rightarrow0\text{ as }m\rightarrow\infty.
\label{cond interTh}%
\end{equation}
By Hopf's ergodic theorem for Markov operators (Theorem 11.4 in Eisner et al.
2015)
\[
f_{m}\rightarrow0\text{ }\pi-a.s.\text{ so }f_{m}g_{f}\rightarrow0\text{
\ }\pi-a.s.\text{ }%
\]
and also, because by condition (\ref{Qnew}), $\ (f_{m}g_{f})_{m\geq1}$ is
uniformly integrable, it follows that
\[
f_{m}g_{f}\rightarrow0\text{ in }\mathbb{L}_{1}(\pi).\text{
\ \ \ \ \ \ \ \ \ \ \ \ \ \ \ \ \ \ \ \ \ \ \ \ \ \ }\square
\]

\bigskip

\textbf{Proof of Corollary \ref{corstrong}. }

Note that, by the triangle inequality, (\ref{strong}) implies
(\ref{cond interTh}) and the proof of Theorem \ref{Th-quenched} applies.

\bigskip

P\textbf{roof of Corollary \ref{cor-quenched}.}

We start from (\ref{cond interTh}) and apply H\"{o}lder's inequality, so%
\[
\mathbb{E}_{\pi}(|f_{m}g_{f}|)\leq\mathbb{E}_{\pi}^{1/p}(|f_{m}|^{p}%
)\mathbb{E}_{\pi}^{1/q}(|g_{f}|^{q}).
\]
By the mean ergodic theorem for the Dunford-Schwartz operators on a Banach
space (see Theorem 8.18 in Eisner et al., 2015) $\mathbb{E}_{\pi}(|f_{m}%
|^{p})\rightarrow0$ as $m\rightarrow\infty,$ and the result follows. Also note
that we can take $p=\infty$ and $q=1.$\ $\square$

\bigskip

\textbf{Proof of Corollary \ref{rev}.}

We shall verify the condition\ of Corollary \ref{cor-quenched}.\textbf{ }If
$f\in(I-Q)\mathbb{L}_{q}(\pi)$ there is $h\in\mathbb{L}_{q}(\pi)$ such that
$f=(I-Q)h.$

Then, by H\"{o}lder's inequality%
\[
\mathbb{E}_{\pi}(\sup_{n}|(I+Q+...+Q^{n-1})f|^{q})=\mathbb{E}_{\pi}(\sup
_{n}|(I-Q^{n})h|^{q})\leq2^{q-1}[\mathbb{E}_{\pi}(|h|^{q})+\mathbb{E}_{\pi
}(\sup_{n}(|Q^{n}h|^{q})],
\]
By the Stein Theorem (see Stein, 1961), $\sup_{n}|Q^{n}h|$ is in
$\mathbb{L}_{q}(\pi)$ and there is a constant $K$ such that $\mathbb{E}_{\pi
}\sup_{n}(|Q^{n}h|^{q})\leq K\mathbb{E}_{\pi}(|h|^{q}).$ Therefore $g_{f}$ is
in $\mathbb{L}_{q}(\pi)$ and we can apply Corollary \ref{cor-quenched} to
obtain the result. $\ \square$

\bigskip

\textbf{Proof of the Remark \ref{example}}.

It is convenient to specify this example in terms of a stationary process
defined by a dynamical system. The proof of this remark follows by analyzing
the example given in Durieu and Voln\'{y} (2008) and Durieu (2009).

We consider an ergodic dynamical system $(\Omega,{\mathcal{A}},\mu,T),$ with
$\mu$ nonatomic and strictly positive entropy. Let $\mathcal{B}$ and
$\mathcal{C}$ be two independent sub-sigma algebras of ${\mathcal{A}}$ such
that $T^{-1}\mathcal{C}=\mathcal{C}.$ Let $(e_{i})_{i\in\mathbb{Z}}$ be a
sequence of independent identically distributed Rademacher random variables
with parameter $1/2$, measurable with respect to $\mathcal{B}$ and denote by
${\mathcal{F}}_{0}$ the $\sigma$-algebra generated by $\mathcal{C}$ and
$(e_{i})_{i\leq0}$. We consider an increasing sequence of integers $(N_{k}),$
and mutually disjoint sets $(A_{k})_{k\in{\mathbb{Z}}}$, $A_{k}\in\mathcal{C}$
such that

(1) $\frac{2}{3}\rho_{k}\leq\mu(A_{k})\leq\rho_{k}\text{ for all }%
k\in{\mathbb{N}}^{\ast}$ where $\rho_{k}=a^{k}$ for $0<a\leq1/4.$

and

(2) for all $k\in{\mathbb{N}}$ and all $i,j\in\{0,\dots,N_{k}\}$, $\mu
(T^{-i}A_{k}\Delta T^{-j}A_{k})\leq\varepsilon_{k}$ where $(\varepsilon_{k})$
will be selected later.

The existence of the sequence $(A_{k})_{k\in{\mathbb{Z}}}$ with the above
properties was explained in Lemma 2 of Durieu and Voln\'{y} (2008).

The function $f$ is then defined as
\begin{equation}
\text{$f=\sum_{k\geq1}e_{-N_{k}}\mathbf{1}_{A_{k}}$}\,. \label{deffunction}%
\end{equation}
The function $f$ defined in (\ref{deffunction}) is centered, ${\mathcal{F}%
}_{0}$-measurable and bounded.

For any $i\in{\mathbb{Z}},$ let now $X_{i}=f\circ T^{i}$ . This sequence is
adapted to the stationary and nondecreasing sequence of $\sigma$-algebras
$({\mathcal{F}}_{i})_{i\in{\mathbb{Z}}}$ where ${\mathcal{F}}_{i}%
=T^{-i}({\mathcal{F}}_{0})$. Note that the sequence $(e_{i})_{i\in\mathbb{Z}}$
is adapted to $({\mathcal{F}}_{i})_{i\in{\mathbb{Z}}}$ and ${\mathbb{E}}%
(e_{i}|{\mathcal{F}}_{0})=e_{i}\mathbf{1}_{i\leq0}$ almost surely. Also, for
all $k$ and $i$, $\mathbf{1}_{A_{k}}\circ T^{i}$ is ${\mathcal{F}}_{0}%
$-measurable and the $e_{i}$'s and the $\mathbf{1}_{A_{k}}$'s are independent.
Clearly, for any $i\in{\mathbb{N}}$,
\begin{align}
{\mathbb{E}}(X_{i}|{\mathcal{F}}_{0})  &  =\sum_{k\geq1}e_{-N_{k}+i}%
\mathbf{1}_{i\leq N_{k}}\mathbf{1}_{T^{-i}(A_{k})}\label{consta1CE}\\
&  =\sum_{k\geq1}e_{-N_{k}+i}\mathbf{1}_{i\leq N_{k}}\mathbf{1}_{A_{k}}%
+\sum_{k\geq1}e_{-N_{k}+i}\mathbf{1}_{i\leq N_{k}}(\mathbf{1}_{T^{-i}%
(A_{k})\backslash A_{k}}-\mathbf{1}_{A_{k}\backslash T^{-i}(A_{k})}).\nonumber
\end{align}
So, by using the fact that the $e_{j}$'s and $f$ are bounded by one, and
selecting $N_{k},\varepsilon_{k}$ such that $\sum_{k\geq1}N_{k}\varepsilon
_{k}<\infty,$ we obtain
\begin{gather}
\sum_{i\geq1}{\mathbb{E}}|\sum_{k\geq1}e_{-N_{k}+i}\mathbf{1}_{i\leq N_{k}%
}(\mathbf{1}_{T^{-i}(A_{k})\backslash A_{k}}-\mathbf{1}_{A_{k}\backslash
T^{-i}(A_{k})})|\leq\sum_{i\geq1}\sum_{k\geq1}\mathbf{1}_{i\leq N_{k}}%
[\mu(T^{-i}(A_{k})\Delta A_{k})]\label{consta2CE}\\
\leq\sum_{k\geq1}N_{k}\varepsilon_{k}<\infty\,.\nonumber
\end{gather}
Therefore, since $f$ is bounded, by (\ref{consta1CE}) and (\ref{consta2CE}%
),\ in order to show that $\sum_{i\geq0}{\mathbb{E}}|f{\mathbb{E}}%
(X_{i}|{\mathcal{F}}_{0})|=\infty$ holds, it is enough to show that
\begin{equation}
\sum_{i\geq1}{\mathbb{E}}|f\sum_{k\geq1}e_{-N_{k}+i}\mathbf{1}_{i\leq N_{k}%
}\mathbf{1}_{A_{k}}|=\infty\,. \label{CEbut1}%
\end{equation}
By the fact that $(A_{k})$ are disjoint $\ $%
\begin{gather*}
\sum_{i\geq1}{\mathbb{E}}|f\sum_{k\geq1}e_{-N_{k}+i}\mathbf{1}_{i\leq N_{k}%
}\mathbf{1}_{A_{k}}|=\sum_{i\geq1}{\mathbb{E}}|\text{$\sum_{u\geq1}e_{-N_{u}%
}\mathbf{1}_{A_{u}}$}\sum_{k\geq1}e_{-N_{k}+i}\mathbf{1}_{i\leq N_{k}%
}\mathbf{1}_{A_{k}}|=\\
\sum_{i\geq1}{\mathbb{E}}|\text{$\sum_{k\geq1}e_{-N_{k}}$}e_{-N_{k}%
+i}\mathbf{1}_{i\leq N_{k}}\mathbf{1}_{A_{k}}|=\sum_{i\geq1}\text{$\sum
_{k\geq1}\mathbb{E}|e_{-N_{k}}$}e_{-N_{k}+i}\mathbf{1}_{i\leq N_{k}}%
\mathbf{1}_{A_{k}}|=\\
\sum_{i\geq1}\sum_{k\geq1}\mathbf{1}_{i\leq N_{k}}\mu(A_{k})\geq\frac{2}%
{3}\sum_{i\geq1}\sum_{k\geq1}\mathbf{1}_{i\leq N_{k}}\rho_{k}=\frac{2}{3}%
\sum_{k\geq1}N_{k}\rho_{k}.
\end{gather*}
On the another hand, by (\ref{consta1CE}) and (\ref{consta2CE}),
\begin{equation}
\mathbb{E}\sup_{n}|\sum\limits_{i=1}^{n}{\mathbb{E}}(X_{i}|{\mathcal{F}}%
_{0})|\leq\sum_{k\geq1}\mathbb{E}\sup_{n}|\sum\limits_{i=1}^{n\wedge N_{k}%
}e_{-N_{k}+i}\mathbf{1}_{A_{k}}|+\sum_{k\geq1}N_{k}\varepsilon_{k}.
\label{rel4}%
\end{equation}
By the fact that $(e_{i})$'s and $(A_{k})$'s are independent and by Doob's
maximal inequality we obtain%
\begin{align*}
\sum_{k\geq1}\mathbb{E}\sup_{n}|\sum\limits_{i=1}^{n\wedge N_{k}}e_{-N_{k}%
+i}|\mathbf{1}_{A_{k}}  &  =\sum_{k\geq1}\mathbb{E}\max_{1\leq j\leq N_{k}%
}|\sum\limits_{i=1}^{j}e_{-N_{k}+i}|\mu(A_{k})\\
&  \leq\sum_{k\geq1}\mathbb{E}\max_{1\leq j\leq N_{k}}|\sum\limits_{i=1}%
^{j}e_{-N_{k}+i}|\rho_{k}\leq\sum_{k\geq1}\sqrt{N_{k}}\rho_{k}.
\end{align*}
To finish the proof of this remark we have to select sequences such that
$\sum_{k\geq1}N_{k}\varepsilon_{k}<\infty,$ $\sum_{k\geq1}N_{k}\rho_{k}%
=\infty$ and $\sum_{k\geq1}\sqrt{N_{k}}\rho_{k}<\infty.$

This selection is possible. For instance, we can take $\rho_{k}=4^{-k}$,
$N_{k}=4^{k}$ and $\varepsilon_{k}=8^{-k}.$

\textbf{Proof of the Remark \ref{strong mixing}. Application to strong mixing
sequences.}

We shall apply now Corollary \ref{corstrong} to strong mixing sequences.

For the random variable $X$, define the "upper tail" quantile function $q$ by
\[
q(u)=\inf\left\{  t\geq0:{\ \mathbb{P}}\left(  |X_{0}|>t\right)  \leq
u\right\}  .
\]

Relevant to this application is the following lemma.

\begin{lemma}
\label{inequality2} Let $(\Omega,\mathcal{A},{\mathbb{P}})$ be a probability
space and $\mathcal{M}$ be a $\sigma$-algebra of $\mathcal{A}$. Let $X$ and
$Y$ be two square integrable identically distributed random variables. Denote
by $q$ their common quantile function. Then%
\[
{\mathbb{E}}|X{\mathbb{E}}(Y|\mathcal{M})|\leq3\int_{0}^{\bar{\alpha}}%
q^{2}du,
\]
where%
\[
\bar{\alpha}=\bar{\alpha}(Y,\mathcal{M})=\sup_{t\in{\mathbb{R}}}%
\mathbb{E}|\mathbb{P}(Y\leq t|\mathcal{M})-\mathbb{P}(Y\leq t)|.
\]

\end{lemma}

Inspired by the proof of Lemma 2 in Merlev\`{e}de et al. (1997), this lemma
can be obtained directly, by truncation arguments. It can also be obtained by
using Lemma 4 in Merlev\`{e}de and Peligrad (2006), combined with Rio's
covariance inequality (Theorem 1.1 in Rio, 2000). The proof is left to the reader.

\qquad Let $(X_{i})_{i\in{\mathbb{Z}}}$ be a stationary sequence of real
valued random variables. We shall interpret it as a function of a Markov chain
$\xi_{k}=(X_{j},j\leq k),$ $f(\xi_{k})=X_{k},$ and define the $\sigma$-algebra
${\mathcal{F}}_{0}=\sigma(X_{i},i\leq0)$. For any $k\in{\mathbb{N}}$ also
define
\[
\bar{\alpha}_{k}=\bar{\alpha}(X_{k},{\mathcal{F}}_{0}).
\]
Recall that the strong mixing coefficient of Rosenblatt (1956), defined by
\[
\alpha_{k}=\sup_{A\in\sigma(Y_{k}),B\in{\mathcal{F}}_{0}}|{\mathbb{P}}(A\cap
B)-{\mathbb{P}}(A){\mathbb{P}}(B)|,
\]
is such that $\bar{\alpha}_{k}\leq2\alpha_{k}.$ (see page 8 in Rio, 2000).

By using our Corollary \ref{corstrong} we shall establish the following result
(see also Corollary 3.5 in Dedecker et al., 2014).

For any nonnegative random variable $Z$, we define the quantile function
$q_{Z}$ of $Z$ by $q(u)=\inf\{t\geq0:{\mathbb{P}}(|Z|>t)\leq u\}$.

\begin{proposition}
Assume $(X_{i})_{i\in{\mathbb{Z}}}$ is a stationary and ergodic sequence of
random variables and $|X_{0}|$ has quantile function $q.$ Also assume
\begin{equation}
\sum_{j\geq1}\int_{0}^{\bar{\alpha}_{j}}q^{2}du<\infty. \label{strongRio}%
\end{equation}
Then the quenched functional CLT\ holds.
\end{proposition}

\textbf{Proof.} Note that $\mathbb{E}_{\pi}|(Q^{m}f)(Q^{j}f)|\leq
\min(\mathbb{E}|f(\xi_{m})(Q^{j}f)(\xi_{0})|,$ $\mathbb{E}|f(\xi_{j}%
)(Q^{m}f)(\xi_{0})|).$ So, by Lemma \ref{inequality2} we obtain%
\begin{equation}
\sum_{j\geq1}\mathbb{E}|(Q^{m}f)(Q^{j}f)|\leq3\sum_{j\geq1}\min(\int
\nolimits_{0}^{\bar{\alpha}_{j}}q^{2}du,\int\nolimits_{0}^{\bar{\alpha}_{m}%
}q^{2}du). \label{3}%
\end{equation}
If we impose condition (\ref{strongRio}), this condition implies $\bar{\alpha
}_{m}\rightarrow0$ and also allows us to apply the discrete Lesbesgue
dominated theorem in (\ref{3}). So condition (\ref{strong}) is satisfied and
the result follows. $\ \square$

We easily recognize condition (\ref{strongRio}) as being the usual condition,
optimal in some sense, used in the context of invariance principles for
strongly mixing sequences (see Doukhan et al., 1994).

Note that $X_{0}$ is distributed as $q(U)$ where $U$ is a uniform random
variable. Therefore we can give sufficient conditions for the validity of
(\ref{strongRio}) in terms of moments of $X_{0}$ and mixing rates.

For instance if $X_{0}$ is almost surely bounded by a constant, condition
(\ref{strongRio}) is satisfied as soon as $\sum_{j\geq1}\bar{\alpha}%
_{j}<\infty.$ If for a $\delta>0$ we have $\mathbb{E(}|X_{0}|^{2+\delta
})<\infty,$ then condition (\ref{strongRio}) is satisfied provided
$\sum_{j\geq1}j^{2/\delta}\bar{\alpha}_{j}<\infty$ (see Doukhan et al.,
1994)$.$

\bigskip

\noindent\textbf{Acknowledgements.} The authors are indebted to Florence
Merlev\`{e}de for helpful discussions. Many thanks are going to the referee
for carefully reading the manuscript and for numerous suggestions which
improved the presentation of the paper. This research was supported in part by
a Charles Phelps Taft Memorial Fund grant and the NSF\ grant DMS-1512936.

\end{document}